\documentclass[12pt]{amsart}

\date{\today}
\usepackage{amsmath,amsthm,amssymb}
\newtheorem{lemma}{Lemma}[section]
\usepackage{mathptmx}
\usepackage{hyperref}
\usepackage{mathrsfs}
\usepackage[T1]{fontenc}
\usepackage{bbm}
\author{Hannah Cairns}
\address{Department of Mathematics and Statistics\\\newline
    \hbox{} \hskip 1em Burnside Hall, Room 1005\\\newline
    \hbox{} \hskip 1em 805 Sherbrooke Street West\\\newline
    \hbox{} \hskip 1em Montreal, Quebec, Canada, H3A 0B9}
\email[H.~Cairns]{hannah.abigail.cairns@gmail.com}

\def\angry#1.{\noindent {\bf #1.}}
\def\1{\mathbbm{1}}
\def\curlyh{{\mathcal{H}}}
\def\cuaver#1{{#1}^{\it ca}}
\title[Smash sum is unique]{The smash sum is the unique sum of open sets satisfying a natural list of axioms}
\newtheorem{corollary}[lemma]{Corollary}
\newtheorem{theorem}[lemma]{Theorem}

\makeatletter
\newtheoremstyle{smalllemma}
  {6pt}
  {6pt}
  {\narrower \small}
  {}
  {\bfseries}
  {.}
  {.5em}
  {}
\makeatother
\theoremstyle{smalllemma}

\def\setdelta{\mathbin{\text{\footnotesize $\Delta$}}}
\def\rad{\mathop{\text{\rm rad\hskip 0.09em}}}
\def\bulk#1{\mathscr{B}(#1)}
\def\sumsupdf{\mathbin{\&}}
\long\def\eat#1{}
\def\setdelta{\mathbin{\text{\footnotesize $\Delta$}}}

\begin{document}
\maketitle
\begin{abstract}
A sum of open sets is a map taking two bounded open sets $A, B$ and producing a new open set $A \oplus B$. We prove that, up to sets of measure zero, there is only one such sum satisfying a natural list of axioms. It is the scaling limit of the Diaconis-Fulton smash sum.
\end{abstract}
\tableofcontents

\begin{section}{Sums of open sets and physical models}
\begin{subsection}{Sums of open sets}\label{introduction}
A \emph{sum of open sets} is a binary operator $\oplus$ that takes two bounded open subsets $A, B \subseteq \mathbb R^d$ and produces a new open set $A \oplus B \subseteq \mathbb R^d$, which may be unbounded.

There are many such sums, so we will add some requirements.
A good sum should be monotone, commutative, and associative.
Let $A, B, C$ be bounded open sets. We want:

\begin{itemize}
\item Monotonicity. $A \subseteq A \oplus B$, and
if $A \subseteq B$, then $A \oplus C \subseteq B \oplus C$.
\item Commutativity. $A \oplus B = B \oplus A$.

\item Associativity when it makes sense. If $A \oplus B$ and $B \oplus C$ are both bounded, then $(A \oplus B) \oplus C = A \oplus (B \oplus C)$.
\end{itemize}

A good sum should also have some symmetry properties. If $A$ is an open set and $x \in \mathbb R^d$, let $A + x = \{a + x: a \in A\}$.
Then we ask for

\begin{itemize}
\item Translation invariance. $(A + x) \oplus (B + x) = (A \oplus B) + x$.
\end{itemize}

Let $\mathcal H \subseteq O(d)$ be the group of isometries that fix the unit cube. Each such isometry is a permutation of coordinates followed possibly by changes of sign in some coordinates. We call these \emph{cubic isometries}, and we ask for the sum to be invariant with respect to these isometries.\footnote{We may use a different symmetry group here. We only need two conditions for the proof to work: first, $\curlyh$ should have no fixed points except $0$, which is enough to get Lemma~\ref{diameter}. Second, the statement of Lemma~\ref{poly} should be true.

For example, the proof works for triangular or tetrahedral symmetry, but not for dihedral symmetry in $\mathbb R^3$, because if the plane of reflection is the $xy$ plane, then the average of $z^2$ under all isometries is $z^2 \ne (x^2+y^2+z^2)/3$.}

\begin{itemize}
\item Rotation invariance. If $U \in \mathcal{H}$, then $(UA) \oplus (UB) = U(A \oplus B)$.
\end{itemize}

Finally, let $\lambda$ be the usual Lebesgue measure on $\mathbb R^d$. Our last request is

\begin{itemize}
\item Conservation of measure. For any bounded open sets $A, B$,
$$\lambda(A \oplus B) = \lambda(A) + \lambda(B).$$
\end{itemize}

There is at least one sum that satisfies these six requirements, called the ``smash sum,'' defined in \cite{levineperes}. We will also call it the ``continuous Diaconis-Fulton smash sum'' to distinguish it from the discrete version in \cite{diaconisfulton}. The main theorem of this paper is that it is the only sum that satisfies the requirements, up to sets of measure zero.

We delay the definition of the sum to Section~\ref{continuum}. The fact that the definition makes sense is not at all trivial; it is a theorem due to Sakai \cite{sakaiobstacle}. For the convenience of the reader, we give a relatively elementary proof of that theorem in a self-contained supplement.

\end{subsection}
\begin{subsection}{Physical models}\label{physical}

This set of requirements is motivated by mathematical models of particle systems. These models have a boundary: some of the space is occupied by particles, while some of it is empty. When the local density is low, the particles stay in one place, but when the density of particles exceeds some threshold, they move outward and enter new areas. The models we consider are invariant under cubic isometries, at least.

One example of such a process is internal diffusion-limited aggregation, a discrete process where particles walk randomly on a lattice until they find an empty vertex, and they stop there; for details, see \cite{id}, \cite{diaconisfulton}, \cite{levineperes}. Another example, this time a continuous process, is Hele-Shaw flow. Here water, bounded by air or some other fluid, is allowed to move slowly between two parallel plates that are very close together. Water is almost incompressible, so the density is roughly constant, and surface tension is negligible. If more water is injected, the boundary moves outward in a predictable way.

Both of these processes are closely related to the smash sum. The scaling limit of internal diffusion-limited aggregation (and other similar models) is the smash sum, as proven in \cite{levineperes}. The set in Hele-Shaw flow at time $t$ is the same as the set obtained by using smash sum to add many small balls with total mass $t$ centered at the injection point; see \cite{laplacian}, especially Section 3.3.

\eat{\footnote{
    If a ball $B_\varepsilon(x)$ is contained in an open set $A$, then adding the ball is the same as adding a Dirac at $x$ of the same mass. Proof: The two weight functions are
        $w =
            \1_A + \1_{B_\varepsilon(x)}$ and $w' = \1_A + \delta_0 \lambda(B_\varepsilon)$. If $Q$ is a quadrature domain for $w'$, and $s$ is integrable and superharmonic on $Q$, then $\int (w' - w) s \,dx = \lambda(B_\varepsilon) (s(x) - {1 \over \lambda(B_\varepsilon)} \int_{B_\varepsilon(x)} s(y) \,dy)$, and the value in parentheses is at least zero because a superharmonic function is at least equal to its average on balls. But this means that
        $\int w's \,dx
            \ge \int ws \,dx
            \ge \int_{Q} s \,dx$, so $Q$ is also a quadrature domain for $w'$. See Section 2 for the definition of quadrature domains.
}}

If one already knows that there is a sum of open sets associated with these models, then one would expect it to conserve mass, as well as being monotone, commutative, translation invariant (since the lattice becomes infinitely fine in the limit), and invariant under cubic isometries because those are symmetries of the lattice. Moreover, the discrete processes are all ``abelian networks'' in the sense of \cite{blev}, which roughly means that the final state does not depend on the order of events. So one would expect that the sum should be associative, that is, independent of the order of addition.

For this reason, it seems likely that the uniqueness theorem in this paper could be used to re-prove the scaling limits in \cite{levineperes} from a different direction, by proving that a scaling limit exists and then showing that it must be a sum of open sets satisfying the six requirements above. However, we do not attempt this here.

We now begin the proof of uniqueness. We will first play around with these requirements and derive some elementary consequences, and then describe a winning strategy for a solitaire game, which we call ``smash game.''

\end{subsection}
\begin{subsection}{Inflations and boundedness}

Let $B_r(x) := \{y \in \mathbb R^d: |y - x| < r\}$ be the open ball of radius $r$ centred at~$x$.
If $E$ is any set, let $d(x, E) := \inf\{|x - y|: y \in E\}$.
Let the \emph{inflation} of an open set $A$ by $\varepsilon > 0$ be $$A^\varepsilon := \{x \in \mathbb R^d: d(x, A) < \varepsilon\} = \bigcup_{a \in A} B_\varepsilon(a) = \bigcup_{|x| < \varepsilon} A + x.$$
It is clearly an open set also.

For any sum of open sets, the inflation of the sum is contained in the sum of the inflation.
To see this, let $x \in \mathbb R^d$, $|x|<\varepsilon$.
Then $A + x \subseteq A^\varepsilon, B + x \subseteq B^\varepsilon$.
By translation invariance and monotonicity, 
$$(A \oplus B) + x = (A + x) \oplus (B + x) \subseteq A^\varepsilon \oplus B^\varepsilon.$$
Taking the union over all the points $|x|<\varepsilon$, we have the promised inclusion
\begin{equation}\label{inflationinclusion}(A \oplus B)^\varepsilon \subseteq A^\varepsilon \oplus B^\varepsilon.\end{equation} This is the \emph{inflation inclusion}.

\begin{subsubsection}{Deflation}
We let the \emph{deflation} of an open set $A$ be $$A^{-\varepsilon} := \{x \in \mathbb R^d: d(x, A^c) > \varepsilon\}.$$
The deflation is also an open set, because $d(x, A^c)$ is continuous in $x$.

We claim that, if $A$ is open, then $$A^{-\varepsilon} = \bigcap_{|x| \le \varepsilon} A + x.$$
If $y$ is in the left set, then $d(y, A^c) > \varepsilon$, so $y - x \in A$ if $|x| \le \varepsilon$,
and that means that $y$ is in the intersection. On the other hand, if $y$ is not in the left set, let
$z$ be a point in the closed set $A^c$ with $d(y,A^c)=d(y,z)\le\varepsilon$.
Then $z \notin A$, so $y \notin A+y-z$ and $y$ is not in the intersection.
This proves the equality.%
\footnote{
That is not necessarily true if $A$ is not open; for example, if $A$ is the closed unit ball, then the deflation of $A$ by $1$ is empty, but $\{x: \overline{B_1(x)} \subseteq A\}$ is the point $\{0\}$.
}

If $|x| \le \varepsilon$, then we have $A+x\supseteq{A}^{-\varepsilon}$, $B+x\supseteq{B}^{-\varepsilon}$. By translation invariance and monotonicity, $(A + x) \oplus (B + x) \supseteq A^{-\varepsilon} \oplus B^{-\varepsilon}$. Taking the intersection of that inclusion over all points $|x| \le \varepsilon$, we find that the deflation of a sum contains the sum of the deflations: $$(A \oplus B)^{-\varepsilon} \supseteq A^{-\varepsilon} \oplus B^{-\varepsilon}.$$

\end{subsubsection}
\end{subsection}
\begin{subsection}{The diameter of the sum}
Let $B_r(x)$ be the open ball of radius $r$ centered at the point $x$. Let $B_r := B_r(0)$.
We prove that the sum of two balls is contained in a ball that's not much larger.
\begin{lemma}\label{diameter} There is a constant $N = N_d$ with $B_r \oplus B_r \subseteq B_{Nr}$.\end{lemma}

\angry Proof.
Let $N = 2/((9/8)^{1/d} - 1)$.
Suppose $B_r \oplus B_r \not\subseteq B_{Nr}$. Let $x$ be a point in $B_r \oplus B_r$ with $|x| > N$.
Then the three points $x, 0, -x$ are all in $B_r \oplus B_r$, by monotonicity and cubic symmetry, so the inflation inclusion (\ref{inflationinclusion}) says
\begin{align*}
(B_r)_{Nr/2}\oplus(B_r)_{Nr/2}
    &\supseteq (B_r \oplus B_r)_{Nr/2}
    \\&\supseteq B_{Nr/2}(x) \cup B_{Nr/2}(0) \cup B_{Nr/2}(-x).
\end{align*}
These balls are all disjoint.

The inflation of a ball $B_r$ by $s$ is $B_{r+s}$, so $(B_r)_{Nr/2} = B_{(N/2 + 1)r}$.  By conservation of mass, the measure of the left side is $2(N/2+1)\lambda(B_r)$. The left side contains the right side, so the measure of the left is at least the measure of all three disjoint balls: $$2\left({N\over2}+1\right)^d\lambda(B_r) \ge 3\left({N\over2}\right)^d \lambda(B_r).$$
Therefore, the ratio of the two sides is at least $1$. However, for our chosen value, $2(N/2 + 1)^d / 3(N/2)^d = \frac23 (1+2/N)^d = 3/4$.
This is a contradiction, so our assumption was wrong, and $B_r \oplus B_r$ is a subset of $B_{Nr}$.
\qed

\medskip \angry Remark. In this lemma, we can replace the cubic isometry group by any group of isometries $\curlyh$ that has no fixed points except the origin. Let $$\inf_{|x| = 1} \sup_{U \in \curlyh} |Ux - x| = c >0.$$ Then we inflate everything by $Nrc/2$ so that the three balls are still disjoint, and the lemma still holds with $N = (1/c) \times 2/((9/8)^{1/d} - 1)$.

\begin{corollary}
Any sum satisfying the requirements is bounded.
\end{corollary}

\angry Proof. Let $A, B$ be bounded open sets. Let $r$ be large enough that $A, B \subseteq B_r$. Then $A \oplus B \subseteq B_{Nr}$, which is bounded.
\qed
\medskip

Now that we know this, we can drop the requirement of boundedness in associativity: $A \oplus B$ is always bounded, so $(A \oplus B) \oplus C = A \oplus (B \oplus C)$ for any sets $A, B, C$.

\end{subsection}
\begin{subsection}{A weaker version of the six requirements}

\smallskip

As before, we say that two sets $A,B$ are \emph{essentially equal} if $\lambda(A \setdelta B) = 0$, and we say that $A$ is \emph{essentially contained} in $B$ if $\lambda(A \setminus B) = 0$.

Let $A, B, C$ be bounded open sets. Then a sum of open sets satisfies the requirements \emph{in the essential sense} if:

\begin{itemize}
\item $A$ is essentially contained in $A \oplus B$.
If $A$ is really contained in $B$, then $A \oplus C$ is essentially contained in $B \oplus C$.
\item $A \oplus B$ is essentially equal to $B \oplus A$.
\item If $A \oplus B$ and $B \oplus C$ are bounded, then $(A \oplus B) \oplus C$ is essentially equal to $A \oplus (B \oplus C)$. 
\item If $x \in \mathbb R^d$, then $(A + x) \oplus (B + x)$ is essentially equal to $(A \oplus B) + x$.
\item If $U \in \mathcal H$, then $(UA) \oplus (UB)$ is essentially equal to $U(A \oplus B)$.
\item The measure of the sum $A \oplus B$ is $\lambda(A) + \lambda(B)$.
\end{itemize}

If a sum obeys the requirements in the essential sense, there is an essentially equal sum that really obeys the requirements.
We will prove that in this section as Lemma~\ref{aseventhaxiom}. First, we introduce the idea of a bulky set.

\end{subsection}
\begin{subsection}{Bulky open sets}
\label{bulkyopensets}

If $A$ is an open subset of $\mathbb R^d$, let $\mathscr{U}(A)$ be the set of open sets which are essentially equal to $A$.
We call this the \emph{open equivalence class} of $A$.

We say that an set is \emph{bulky} if it is an open set that contains every other set in its open equivalence class.

\begin{lemma} If $A$ is an open set, then there is exactly one bulky set in $\mathscr{U}(A)$.\label{abulkyset}
\end{lemma}

\angry Proof. First we prove that there is at least one. Let $U$ be the union of all open balls with rational centers and radii that are essentially contained in $A$.
This is a countable union, so $U$ is also essentially contained in $A$.

Let $V$ be any set in $\mathscr{U}(A)$.
Fix $x \in V$.
Let $B$ be a rational open ball with $x \in B \subseteq V$. Then $\lambda(B \setminus A) \le \lambda(V \setminus A) = 0$, so the ball $B$ is essentially contained in $A$. Therefore, $B$ is one of the balls in the union defining $U$,
and it follows that $x \in B \subseteq U$.
Therefore, $U$ contains every set in $\mathscr{U}(A)$.

In particular, it contains $A$. On the other hand, we chose $U$ so that it is essentially contained in $A$. which means that they are essentially equal and share the same equivalence class $\mathscr{U}(A) = \mathscr{U}(U)$.

We have already seen that $U$ contains every set in the open equivalence class of $A$, so it is bulky.
If $\mathscr{U}(A)$ had two bulky sets, they would have to contain each other, which is absurd. Therefore there is exactly one.
\qed

\medskip

Let the unique bulky set that is essentially equal to $A$ be denoted by~$\mathscr{B}(A)$.
Two open sets $A$, $B$ are essentially equal if and only if $\bulk A = \bulk B$.

\begin{lemma} \label{really} $A$ is essentially contained in $B$ if and only if $\bulk A \subseteq \bulk B$.
\end{lemma}

\angry Proof. \emph{If:}
If $\bulk{A} \subseteq \bulk{B}$, then $A \subseteq \bulk{B}$, and $\bulk{B}$ is essentially equal to $B$, so $A$ is essentially contained in $B.$

\emph{Only if:} If $A$ is essentially contained in $B$, then $\bulk{A} \cup B$ is essentially equal to $B$.
By the definition of a bulky set, $\bulk{A} \cup B$ is contained in $\bulk{B}$, and so certainly $\bulk{A} \subseteq \bulk{B}$.
\qed

\medskip
Here are some other easy consequences which we will use later.
\begin{itemize}
\item The measures of the sets $A$ and $\bulk{A}$ are the same.
\item If $x \in \mathbb R^d$, then $\bulk{A + x} = \bulk A + x$.
\item If $U \in \mathcal H$, then $\bulk{UA} = U \bulk A$.
\item If $A$ is bounded, then $\bulk A$ is also bounded.
\item $\bulk A$ is contained in the topological closure of $A$.
\end{itemize}

Some remarks to clarify the idea: A point is in $\bulk A$ if and only if there is a ball containing that point that's essentially contained in $A$. We can get examples of non-bulky sets by taking open sets and subtracting closed sets of measure zero. For example, $(0, 1) \setminus \{\frac12\}$ isn't bulky.\footnote{But there are some non-bulky sets that aren't open sets minus closed sets of measure zero. Let $A$ be an open set which is dense in the square $[0, 1]^2$ but has measure only $1/2$. Let $F_n := \{x \in A: d(x, A^c) \le 1/n\} \cap \{x/2^n: x \in \mathbb Z^2\}$.
Then every point in $A$ has a neighbourhood that intersects only finitely many $F_n$, so
$A \setminus \bigcup F_n$ is open and essentially equal to $A$.
If $A \setminus \bigcup F_n = E \setminus C$ for some bulky open set $E$ and closed set $C$, then $C$ has to contain all the points in the closed sets $F_n$, so it has to contain the boundary of $A$, which has measure $1/2$.}

We'll need this lemma in the next section:
\begin{lemma}
    \label{bracketsinside}
    If a sum of open sets obeys the requirements in the essential sense, then
    $\bulk{A \oplus B} = \bulk{\bulk{A} \oplus B} = \bulk{\bulk{A} \oplus \bulk{B}}$.
\end{lemma}

\angry Proof.
We have $A \subseteq \bulk{A}$, so by the requirement of essential monotonicity, $A \oplus B$ is essentially contained in ${\bulk{A} \oplus B}$.

Bulking does not add measure, so both sets
have measure $\lambda(A) + \lambda(B)$.
Therefore, $A \oplus B$ and $\bulk{A} \oplus B$ are essentially equal.
By earlier remarks, the bulkings are really equal:
$\bulk{A \oplus B} = \bulk{\bulk A \oplus B}$.
The same proof works on the right-hand side to give the second inequality.
\qed

\begin{subsubsection}{How to modify a sum to restore the strong requirements}
\label{modifysum}

Suppose that we have a sum of open sets $\oplus$ that obeys the requirements in the essential sense.
Let the \emph{bulking} of the sum be the map $A, B \mapsto \bulk{A \oplus B}$.
This is a new sum of open sets, and it satisfies the six requirements, by the lemma:

\begin{lemma} \label{aseventhaxiom}If a sum of open sets $\oplus$ satisfies the six requirements in the essential sense, then the bulking of that sum really satisfies the six requirements, and is essentially equal to the original sum.
\end{lemma}

\angry Note. Here we say that two sums $\oplus, \oplus'$ are \emph{essentially equal} if, for every open sets $A, A', B, B'$ with $\lambda(A \setdelta A') = \lambda(B \setdelta B') = 0$, the two sums $A \oplus B$ and $A' \oplus' B'$ are essentially equal.

\medskip
\angry Proof.
The proof that the requirements of monotonicity, commutativity, translation and rotation invariance, and conservation of mass are satisfied is merely to put $\bulk{\cdot}$ around both sides of each essential inclusion or equality, and then use Lemma~\ref{really}.

The associativity property causes a little trouble.
Our assumption is that the original sum is essentially associative:
\begin{itemize}
\item Essential associativity when bounded. If $A, B, C$ are bounded open sets and $A \oplus B, B \oplus C$ are bounded, then $(A \oplus B) \oplus C$ is essentially equal to $A \oplus (B \oplus C)$.
\end{itemize}
(We have not yet proved that an \emph{essential} sum is always bounded!)
We have to prove that the bulking of the sum is associative:
\begin{itemize}
\item Associativity for the bulky sum. If $A, B, C$ are bounded open sets and $\bulk{A \oplus B}, \bulk{B \oplus C}$ are bounded, then $\bulk{\bulk{A \oplus B} \oplus C} = \bulk{A \oplus \bulk{B \oplus C}}$.
\end{itemize}

If $\bulk{A \oplus B}$ is bounded, then certainly $A \oplus B$ is also bounded, and similarly for $B \oplus C$.
Therefore, we can use essential associativity, and $(A \oplus B) \oplus C$ is essentially equal to $A \oplus (B \oplus C)$.

So, $\bulk{(A \oplus B) \oplus C} = \bulk{A \oplus (B \oplus C)}$ for any three open sets $A, B, C$ with $A \oplus B$, $B \oplus C$ bounded.

\medskip
Once we know this lemma, we use Lemma~\ref{bracketsinside}, Lemma~\ref{really}, and then Lemma \ref{bracketsinside} again to get
\begin{align*}
\bulk{\bulk{A \oplus B} \oplus C} &= \bulk{(A \oplus B) \oplus C}
\\&= \bulk{A \oplus (B \oplus C)}
\\&= \bulk{A \oplus \bulk{B \oplus C}}
\end{align*}

Therefore, the bulking of the sum is associative, and the proof for all the other six requirements is straightforward.

Let $A, B, A', B'$ be bounded open sets with $A$ essentially equal to $A'$ and $B$ essentially equal to $B'$. Then Lemma~\ref{bracketsinside} tells us that
$$\bulk{A \oplus B} = \bulk{\bulk{A} \oplus \bulk{B}} = \bulk{\bulk{A'} \oplus \bulk{B'}} = \bulk{A' \oplus B'},$$
so the bulked sum $\bulk{A \oplus B}$ is essentially equal to $A' \oplus B'$.
\qed

\medskip

\end{subsubsection}
\begin{subsubsection}{From now on, we assume our sum is bulky}

From now on, we will assume that we have made the replacement described in Section~\ref{modifysum}, and our new sum of open sets satisfies all six requirements, plus a seventh:
\begin{itemize}
\item Bulkiness. The sum of any two sets $A \oplus B$ is bulky.
\end{itemize}
We'll prove that there is only one sum that satisfies all seven requirements, namely the smash sum.

If we have a sum $\oplus'$ of open sets that satisfies the six requirements in the essential sense, the sum is \emph{essentially equal to} the smash sum, in the sense that $A \oplus B$ is essentially equal to $A \oplus' B$ for any bounded open sets $A, B$.

\newpage

\end{subsubsection}
\begin{subsubsection}{Properties of bulky sums}
In this subsection, we assume that $\oplus$ is a bulky sum that satisfies the six requirements in the strong sense, and we prove some lemmas that will be useful later.

\begin{lemma}\label{disj}
If the open sets $b_1, \ldots, b_n$ are disjoint and bounded, then the sum $b_1 \oplus \cdots \oplus b_n$ is the bulking of the union, $\bulk{b_1 \cup \cdots \cup b_n}$.
\end{lemma}

\angry Proof. The sum $A \oplus B$ must contain both $A$ and $B$ by monotonicity, and its mass must be $\lambda(A) + \lambda(B)$, so it is essentially equal to $A \cup B$.
By Lemma~\ref{abulkyset}, it must be $\bulk{A \cup B}$.\qed

\begin{corollary}\label{disjess}
The sum $b_1 \oplus \cdots \oplus b_n$ is essentially equal to $b_1 \cup \cdots \cup b_n$.
\end{corollary}

\angry Proof. By definition, the bulking of a set is essentially equal to that set, so the bulking of the union is essentially equal to the union. \qed

\begin{lemma}\label{boundy}
Let $E$ be a bounded open set. Let $C$ be bounded and open.
If $\partial C$ has measure zero, then $$(E \setminus \overline{C}) \oplus (E \cap C) = \bulk{E}.$$
\end{lemma}

\angry Proof. Let $E' = (E \setminus \overline{C}) \cup (E \cap C)$.
By Lemma~\ref{disj}, $(E \setminus \overline{C}) \oplus (E \cap C) = \bulk{E'}.$
The set difference between $E'$ and $E$ is contained in $\partial C$, which has measure zero, so $E'$ is essentially equal to $E$ and $\bulk{E'} = \bulk{E}$.
\qed

\end{subsubsection}
\end{subsection}
\end{section}
\begin{section}{The definition of the Diaconis-Fulton sum}\label{continuum}

\begin{subsection}{Preliminary setup}

Let $\Omega \subseteq \mathbb R^d$ be an open set and $Q \subseteq \Omega$ be an open subset.

Recall that a function $s$ taking values in $\mathbb R \cup \{-\infty\}$ is \emph{superharmonic} on~$Q$ if it is locally integrable and lower semicontinuous on $Q$, and
$${1 \over \lambda(B_r)} \int_{B_r(x)} s(y) \,dy \le s(x)$$
for $x \in Q$ and sufficiently small $r > 0$.

\begin{subsubsection}{Quadrature domains}

Let $w \ge 0$ be a bounded measurable function.
A {\it quadrature domain} for $w$ is an open set $Q$
which essentially contains the set $\{w > 0\}$, and which has the property that
$$\int_Q s\,dx \le \int sw\,dx$$
for every function $s$ that is both superharmonic and integrable on $Q$.

This property has an intuitive interpretation that is related to the physical description in Section~\ref{physical}. First of all, to set things up, the basic theory of superharmonic functions says that the average on any ball $B_r(x) \subset Q$ is no larger than $s(x)$. Let $\varphi_r(x) = \1_{|x| < r} / \lambda(B_r)$, then $s(x) \ge \int s(y) \varphi_r(x - y) \,dy$ as long as $r \le d(x, Q^c)$.
If $\rho(x)$ is a measurable function on $Q$ with $0 < \rho(x) < d(x, Q^c)$, then $w'(x) = \int w(y) \varphi_{\rho(y)} \,dy$ is a new weight function where the weight in $w$ has been spread outward in a radially symmetric way.

By Fubini's theorem, the integral against this spread-out weight is no larger than the integral against the original weight:
$$\int s(x)w'(x) \,dx = \iint s(x) \varphi_{\rho(x-y)} w(y) \,dy \le \int s(x)w(x) \,dy.$$
This inequality is also true if the weight spreads out in several steps, as long as the disks stay within $Q$.
Now, what we are doing is using this type of inequality as a definition. When we say that $\int s \1_Q \,dx \le \int sw \,dx$ for any function that is superharmonic in $Q$, what we mean by that is: the weight in~$w$ can spread out radially through the set $Q$ so that the new weight is $\1_Q$.
This matches our picture of particles which spread out radially until they reach a certain fixed density.

\end{subsubsection}
\begin{subsubsection}{Existence and uniqueness}

We use theorems from Sakai~\cite{sakaiobstacle}, which are proven in the supplement.

\begin{theorem}[Sakai]\label{sakaisx} If $w, w' \ge 0$ are two bounded measurable weight functions with $w \le w'$, and $Q$ and $Q'$ are quadrature domains for $w$ and $w'$ respectively, then $Q$ is essentially contained in $Q'$.

If $Q$ and $Q'$ are quadrature domains for the same bounded measurable weight function $w \ge 0$, then $Q$ is essentially equal to $Q'$.\end{theorem}

\angry Proof. The first statement is Lemma~3 in the appended supplement, and the second one is Corollary~4 (or follows immediately).

\begin{theorem}[Sakai]\label{smashexists} If $w \ge 0$ is a bounded measurable weight function that is greater than or equal to one on some bounded open set and zero outside it, then there is a bounded quadrature domain for $w$.\end{theorem}

\angry Proof. This is Theorem 33 in the supplement.

\medskip
In particular, if we choose $w = \1_A + \1_B$, where $A, B$ are bounded open sets, then there is a bounded quadrature domain $Q$ for $w$,
and any other quadrature domain for $w$ is essentially equal to $Q$.

\end{subsubsection}
\end{subsection}
\begin{subsection}{Definition of the sum}\label{definingsum}
If $Q$ is a quadrature domain for $w$, then $\bulk{Q}$ is also a quadrature domain for $w$, because every integrable superharmonic function on $\bulk{Q}$ is also integrable and superharmonic on $Q$.

Therefore, every weight function that satisfies the conditions in Theorem \ref{smashexists} has a bulky quadrature domain, which is unique by Theorem \ref{sakaisx} and Lemma~\ref{abulkyset}.

\medskip
\angry Definition. If $A$ and $B$ are bounded open sets, then the Diaconis-Fulton smash sum is the unique bulky quadrature domain for $\1_A + \1_B$.

\medskip
We will denote the smash sum by $A \sumsupdf B$.

\begin{theorem}Diaconis-Fulton smash sum satisfies all the requirements.\end{theorem}

\angry Proof.
Let $A, B$ be bounded open sets. Let $x \in \mathbb R^n$.
If $s$ is an integrable superharmonic function on $(A \sumsupdf B) + x$, then
$$\int_{(A \sumsupdf B)+x} s \,dy \le \int_{A+x} s \,dy + \int_{B+x} s\,dy,$$
so $(A \sumsupdf B) + x$ is a bulky quadrature domain for $\1_{A+x} + \1_{B+x}$. By uniqueness, it's equal to $(A + x) \sumsupdf (B + x)$.
Therefore the sum is translation invariant.
Rotation and reflection invariance follows in the same way, and so does commutativity.

Conservation of mass is easy: $\pm1$ is harmonic, so
$$\int_{A \sumsupdf B} 1 \,dx \le \int_A 1 \,dx + \int_B 1 \,dx \quad\text{and}\quad \int_{A \sumsupdf B} -1 \,dx \le \int_A -1 \,dx + \int_B -1 \,dx.$$
Therefore, $\lambda(A \sumsupdf B) = \lambda(A) + \lambda(B)$ and the sum conserves mass.

Let $A, B, C$ be bounded open sets. If $s$ is an integrable superharmonic function on $(A \sumsupdf B) \sumsupdf C$, then
$$\int_{(A \sumsupdf B) \sumsupdf C} s \,dx \le \int_{A \sumsupdf B} s \,dx + \int_C s \,dx \le \int_A s \,dx + \int_B s \,dx + \int_C s \,dx.$$
Therefore, $(A \sumsupdf B) \sumsupdf C)$ is a quadrature domain for $\1_A + \1_B + \1_C$,
but so is $A \sumsupdf (B \sumsupdf C)$.
So they're equal and the sum is associative.

If $w \le w'$ are two weight functions that satisfy the conditions in Theorem \ref{smashexists}, then by Theorem \ref{sakaisx}, the bulky quadrature domain of $w$ is essentially contained in the bulky quadrature domain of $w'$. By Lemma \ref{really}, it is really contained.

Let $A, B$ be bounded open sets. Then $A, A \sumsupdf B$ are bulky quadrature domain for $\1_A \le \1_A + \1_B$, so $A \subseteq A \sumsupdf B$. Let $A, B, C$ be bounded open sets with $A \subseteq C$. Then $A \sumsupdf B, C \sumsupdf B$ are bulky quadrature domains for $\1_A + \1_B \le \1_C + \1_B$, so $A \sumsupdf B \subseteq C \sumsupdf B$. So monotonicity holds, and that's the last one.
\qed

\medskip

In the rest of this paper, we will prove that any sum that satisfies all the requirements is essentially equal to the Diaconis-Fulton smash sum.

\end{subsection}
\end{section}
\begin{section}{Uniqueness of the smash sum}

Let $A \oplus B$ be some sum that satisfies the six requirements in the strong sense, and is also bulky.
We will prove that the sum satisfies
\begin{equation}
\label{inequalitytosatisfy}
\int_{A \oplus B} s\,dx \le \int_A s\,dx + \int_B s\,dx
\end{equation}
for any integrable superharmonic function $s$ on $A \oplus B$.

\begin{subsection}{Smash game}
We introduce a solitaire game, {\it smash game}.
Imagine $\mathbb R^d$ is a large dining room table.
A bounded open set $A$ is on the table, and you are holding one bounded open set~$B$ in your hand.

You are given a nonnegative smooth superharmonic function $s$ which is defined on the inflation $(A \oplus B)^{\delta}$ for some small $\delta > 0$. We assume that the derivatives of all orders are bounded. Finally, you are given a small $\varepsilon > 0$, which is the difficulty of the game: it will be harder for smaller $\varepsilon$.

You can make four kinds of moves, which are described in Section \ref{thesmashmoves}.
Your progress in the game is tracked as follows:
\begin{itemize}
\item The {\it current sum} is the sum of the table set and all the hand sets.
This starts out at $A \oplus B$, and all of the moves will decrease it or leave it unchanged.
\item The {\it current mass} is the measure of the current sum, and the
{\it mass in the hand} is the sum of the measures of the hand sets. The starting value of the current mass is $\lambda(A) + \lambda(B) = \lambda(A \oplus B)$.
\item If $A$ is the current table set and $B_1, \ldots, B_m$ are the current hand sets, then the {\it total $s$ integral} is $\int_A s \,dx + \sum_{j=1}^m \int_{B_j} s\,dx$.

The starting value of the total $s$ integral is $\int_A s \,dx + \int _B s \,dx$, which is the right-hand side of the desired inequality (\ref{inequalitytosatisfy}).
\end{itemize}

You lose the game if the current mass is lower than $\lambda(A) + \lambda(B) - \varepsilon$, or if the total $s$ integral is greater than $\int_A s \,dx + \int_B s \,dx + \varepsilon$. You win the game if you have not yet lost, and the mass in your hand is less than $\varepsilon$.

\end{subsection}
\begin{subsection}{The consequence of winning}
If you can win smash game at every difficulty $\varepsilon > 0$, then you get the inequality (\ref{inequalitytosatisfy}) at the start of this section.
\begin{theorem}\label{canwinit}
Let $s$ be a smooth nonnegative superharmonic function on $(A \oplus B)^\delta$ for some $\delta > 0$.
If you can win smash game for any $\varepsilon > 0$, then $$\int_{A \oplus B} s \,dx \le \int_A s \,dx + \int_B s \,dx.$$
\end{theorem}

\angry Proof. Play smash game until you win. Let the table set at the end of the game be $A'$. The current sum decreases monotonically over the course of the game, so the final table set is contained in $A \oplus B$.

The current mass is at least $\lambda(A) + \lambda(B) - \varepsilon$ and the total mass of the hand sets is less than $\varepsilon$,
so the final table mass is at least $\lambda(A) + \lambda(B) - 2\varepsilon$.
Because $s$ is bounded, we get $\int_{A \oplus B} s \,dx \le \int_{A'} s \,dx + 2 \varepsilon \sup s$.
You won the game, so the total $s$ integral isn't more than $\int_As\,dx+\int_Bs\,dx+\varepsilon$.
Therefore
\begin{align*}
\int_{A \oplus B} s\,dx
    &\le \int_{A'} s\,dx + 2 \varepsilon \sup s
    \\&\le \int_A s \,dx + \int_B s \,dx + \varepsilon + 2 \varepsilon \sup s.
\end{align*}
Now let $\varepsilon \to 0$ to get the inequality $\int_{A\oplus{B}}s\,dx\le\int_As\,dx+\int_Bs\,dx$.
\qed

\bigskip
We will now finish the definition of smash game and check that all the moves decrease the current sum, in Section \ref{thesmashmoves}. Then we'll work out a strategy in Sections \ref{startofstrategizing}--\ref{endofstrategizing} that will allow us to win with any function $s$ that meets the requirements of the game. Finally we'll extend the result to any integrable superharmonic function $s$ in Theorem~\ref{superharmonicultimate}, and that will mean that $A \oplus B$ is a quadrature domain.

\end{subsection}
\begin{subsection}{The four moves of the smash game}

\label{thesmashmoves}
Here are the four moves of smash game.
\begin{subsubsection}{Replace a hand set by finitely many disjoint balls}

The first move lets us discard a hand set $B$ and replace it by disjoint open balls $b_1, \ldots, b_n$ so that the union $b_1 \cup \cdots \cup b_n$ is contained in a compact subset of $B$.

The bulking of a set is contained in the closure of the set, as mentioned in Section \ref{bulkyopensets}, so $b_1 \oplus \cdots \oplus b_n \subseteq B$,
and by the monotonicity axiom, the replacement of $B$ by $b_1, \ldots, b_n$ does not increase the current sum.

The total $s$ integral doesn't increase, because the balls are disjoint and their union is contained in $B$, so $\sum_{j=1}^n \int_{b_j} s \,dx = \int_{b_1 \cup \cdots \cup b_j} s \,dx \le \int_B s \,dx$.
The current mass decreases by $\lambda(B) - \lambda(b_1 \cup \cdots \cup b_n)$.

\end{subsubsection}
\begin{subsubsection}{Shrink the table set}

The second move lets us replace the table set $A$ by an open subset $A' \subseteq A$.
By the monotonicity axiom, the current sum and total $s$ integral decrease or stay the same. We may lose some total mass.

\end{subsubsection}
\begin{subsubsection}{Smash part of the table set into the hand set}
We can only use the third move if the table set is bulky.

Let $A$ be the table set. Let $C$ be an open subset of $A$ with a boundary of measure zero. Let~$B$ be a set in the hand.
The third move lets us discard $B$ from the hand and replace it by $B' := (B \oplus C) \setminus \overline{C}$.

Intuitively, we are pressing the hand set $B$ down with a metal stamp in the shape of $C$. We get a new set $B'$ which has the same mass, but that mass is now spread out, and there is a hole shaped like $C$.

Putting $E = A$ in Lemma~\ref{boundy}, we get the equality $(A \setminus \overline{C}) \oplus C = A$. Putting $E = B \oplus C$ in the same lemma gives $B' \oplus C = B \oplus C$. Therefore, $$A \oplus B = (A \setminus \overline{C}) \oplus C \oplus B = (A \setminus \overline{C}) \oplus C \oplus B' =A \oplus B',$$ and the current sum doesn't change. The current mass is the mass of the current sum, so it also does not change.

On the other hand, we don't have much information about the action of the sum, so replacing $B$ by $B'$ will change the total $s$ integral by an unpredictable amount. In the next section, we'll see how to use symmetry to gain some control over the change in the $s$ integral.
\end{subsubsection}
\begin{subsubsection}{Move part of a hand set to the table}
We can only use the fourth move if the boundary of the table set $\partial A$ has measure zero.

Let $A$ be the table set. Let $B$ be a hand set.
The fourth move lets us replace $B$ by $B' = A \cap B$ and change the table set to $A' = A \oplus (B \setminus \overline{A})$.

This move decreases the mass in hand by $\lambda(B) - \lambda(B') = \lambda(B \setminus A)$, so we make progress toward winning as long as $B \setminus A$ has positive measure.

By Corollary \ref{disjess}, the sum of disjoint sets is essentially equal to their union, so $A \oplus (B \setminus \overline{A})$ is essentially equal to $A \cup (B \setminus \overline{A}) = (A \cup B) \setminus \partial A$. We assumed that the boundary $\partial A$ has measure zero, so
$$\int_{A'} s \,dx + \int_{B'} s \,dx = \int_{A \cup B} s \,dx + \int_{A \cap B} s \,dx.$$
and the total $s$ integral doesn't change when we replace $A, B$ by $A', B'$.

The current sum also doesn't change, because
\begin{align*}
    A' \oplus B' &= A \oplus (B \setminus \overline{A}) \oplus B' \\&= A \oplus \bulk{B}
    \\&=A \oplus B.
\end{align*}
The first step is the definition of $A'$, the second step is Lemma \ref{boundy} with $E = B$ and $C = A$, and
the last step comes from Lemma~\ref{bracketsinside}.

To summarize, this move never brings us closer to losing, because the total $s$ integral and the current mass stay the same. We can use it to reduce the mass in hand if we have a hand set which is outside of the table set.

This move is the only one that can get mass out of the hand without decreasing the current mass, so it's necessary to win.

\end{subsubsection}
\end{subsection}
\begin{subsection}{The cookie-cutter smash}
\label{startofstrategizing}

\begin{subsubsection}{A moment to consider our strategy}

How can we win smash game?
We need to rearrange the mass in the hand so that it's outside the table set, and then use the fourth move to get rid of it.

If you didn't have to worry about the total $s$ integral, you could get rid of all the mass in the hand in two moves. Use the third move to smash $B$ into the whole table set, which replaces $B$ by $(B \oplus A) \setminus \overline{A}$. Then use the fourth move to put it all down on the table.

We want to bound the increase in the total $s$ integral, but we also want to move mass outside of the table set. The compromise between these goals is the \emph{cookie-cutter smash}, which is defined below.

\end{subsubsection}
\begin{subsubsection}{Definition of the cookie-cutter set}

Recall that a {\it cubic isometry} is an isometry that preserves the cube $[-1, 1]^d$, and the group of those isometries is called $\mathcal{H}$.
Let $U$ be an element of $\mathcal{H}$. If $x$ is a point in $\mathbb R^d$, let $U_x$ be the map that takes $y \in \mathbb R^d$ to $U(y -x) + x$.

Given a table set $A$, we define a symmetrized version of the set. Let $A$ be an open set. Let $x \in \mathbb R^d, R > 0$.
Then the \emph{cookie-cutter set} for $x, R, A$ is
$$\mathscr C(x, R, A) := B_R(x) \cap \bigcap_{U \in H} U_x A.$$
The intuitive picture in two dimensions is that we start with a disc of clay $B_R(x)$  and then cut out a shape by pressing the cookie-cutter $A$ down in all $|\mathcal H| = 8$ different orientations.

We say that a set $E$ has  \emph{cubic symmetry} around a point $x \in \mathbb R^d$ if $U_x E = E$ for $U \in \mathcal H$.
The cookie-cutter set always has cubic symmetry.

The set $\mathscr C(x, R, A)$ is an open set contained in $A$. If the topological boundary of $A$ has measure zero, then $\partial\mathscr{C}(x,R,A)\subseteq\partial{B_R(x)}\cup\bigcup{U_x\partial{A}}$ also has measure zero.
That means that we can use it to do the third move.

\end{subsubsection}
\begin{subsubsection}{The cookie-cutter smash}

We pick a ball $B_r(x)$ in the hand and a radius $R \in (r, \delta/2N)$, where $N$ is from Lemma \ref{diameter}.
Then we use the third move to smash $B_r(x)$ into $\mathscr C(x, R, A)$. This is a \emph{cookie-cutter smash}, and
the result is that the ball $B_r(x)$ is replaced by the \emph{smash set}
$$E = (B_r(x) \oplus \mathscr C) \setminus \overline{\mathscr C}.$$

We bound the increase in the total $s$ integral with the following lemma.
\begin{lemma}\label{x}
Suppose the boundary of the table set $A$ has measure zero.
Let $B_r(x)$ be a ball in the hand. Let $R \in (r, \delta/2N)$ and $x \in A$.
Let $\mathscr C := \mathscr C(x, R, A)$ be the cookie-cutter set, and
let $E$ be the smash set as above. Then $$\int_Es\,dy \le \int_{B_r(x)} s \,dy+ C_sR^3 \lambda(B_r).$$
Here $C_s$ depends only on $s$.

\end{lemma}

\angry Proof. By translation invariance and Lemma~\ref{diameter}, the smash set is contained in the ball $B_{NR}(x) \subseteq (A \oplus B)^{\delta /2}$,
and $s$ and all its derivatives are uniformly bounded on that set.
Let $x = 0$ for convenience of notation.
We can expand the superharmonic function in a Taylor series around $x = 0$:
$$s(y) = P(y) + Q(y),$$
where $P(y)$ is a polynomial of degree at most two and
$|Q(y)| \le c_s |y|^3$ as long as $|y| < R$. The constant in the bound depends only on the function $s$, because $s$ has bounded derivatives of all orders on the inflated set $A^{\delta / 2} \supseteq B_R$.

If $f$ is any function and $E$ is a set with cubic symmetry,
\begin{align*}
\int_E f\,dy= {1 \over |\mathcal H|} \sum_{U \in \curlyh} \int_{UE} f \,dy
&= \int_E{1 \over |\curlyh|}\sum_{U \in \curlyh} f(Uy) \,dy = \int_E \cuaver f\,dy,
\end{align*}
where $\cuaver f:= |\curlyh|^{-1} \sum_{U \in \curlyh} f(Uy)$. Call $\cuaver f$ the \emph{cubic average of $f$}.
Both of the sets $E$ and $B_r$ have cubic symmetry, so
\begin{align*}
\int_Es\,dy-\int_{B_r}s\,dy&=\int_E\cuaver s\,dy-\int_{B_r}\cuaver s\,dy
\\&=\int_E\cuaver P\,dy-\int_{B_r}\cuaver P\,dy+\int_E\cuaver Q\,dy-\int_{B_r} \cuaver Q\,dy.\end{align*}

Let $a := s(0)$ and $b := -\nabla^2s(0)/2d$. Here $b \ge 0$ because $s$ is superharmonic. The functions $s$ and $P$ have the same derivatives at zero up to second order, so $P(0) = a$ and $-\nabla^2 P(0)/2d = b$. By Lemma \ref{poly} below, the cubic average of $P(y)$ is $\cuaver P = a - b|y|^2$. If $y \notin B_r$, then $a - b|y|^2 \le a - br^2$ and
$$\int_{E \setminus B_r} \cuaver P\,dy \le (a-br^2) \lambda(E \setminus B_r) = (a-br^2) \lambda(B_r \setminus E) \le \int_{B_r \setminus E} \cuaver P\,dy.$$
The middle step uses $\lambda(E) = \lambda(B_r)$, so $\lambda(E \setminus B_r) = \lambda(B_r \setminus E)$.
The last step is like the first: if $y \in B_r$ then $a - br^2 \le a-b|y|^2$.

Adding $\int_{E\cap{B_r}}\cuaver{P}\,dy$ to both sides, we find that $\int_E\cuaver{P}\,dy-\int_{B_r}\cuaver{P}\,dy\le0.$
Therefore, the change of the $s$ integral is bounded above by the difference of the remainder integrals:
\begin{align*}
    \int_E s \,dy - \int_{B_r} s \,dy
        &\le \int_E \cuaver Q \,dy - \int_{B_r} \cuaver Q \,dy
        \\&\le \int_E |\cuaver{Q}| \,dy + \int_{B_r} |\cuaver{Q}| \,dy
        \\&\le 2c_s(NR)^3 \lambda(B_r).
\end{align*}
The last inequality is because $B_r, E$ both have mass $\lambda(B_r)$ and diameter at most $NR$, and $|Q(y)| \le c_s |y|^3$.

So the claimed bound is true with $C_s:=2N^3c_s$.  \qed

\begin{corollary}
If $B_r(x)$ is a ball in the hand and $R \in (r, \delta/2N)$, the cookie-cutter smash increases the $s$ integral by at most $C_sR^3\lambda(B_r)$.
\end{corollary}

\angry Proof. When we do the cookie-cutter smash, it replaces $B_r(x)$ by $E$, and so the $s$ integral changes by $$\int_Es\,dx-\int_{B_r(x)}s\,dx.$$ The lemma tells us that this is bounded above by $C_s \lambda(B_r) R^3$, so the total $s$ integral doesn't increase by more than that.
\qed

\medskip
So the cookie-cutter smash will only increase the $s$ integral by the mass of the ball times the third power of $R$.
We still owe an easy technical lemma:
\begin{lemma}
\label{poly}If $P$ is a polynomial of degree two or less, then $$\cuaver P = P(0) + \nabla^2 P(0) {|y|^2\over{2d}}.$$
\end{lemma}

\angry Proof. Recall the definition of the cubic average, $\cuaver f := |\curlyh|^{-1} \sum_{U\in\curlyh}f(Uy)$.

There are four kinds of monomials of degree less than three: $1$, $y_i$, $y_i^2$, and $y_iy_j$. Here $i, j$ denote distinct indices.

The cubic averages of $y_i$ and $y_i y_j$ are always zero.
Let $V$ be the isometry that takes $(\ldots,y_i,\ldots,y_j,\ldots)$ to $(\ldots,y_j,\ldots,-y_i,\ldots)$.
Let $h$ be one of the above monomials. Then
$$\begin{array}{c@{\hskip 1em}|@{\hskip 1em}c@{\hskip 1em}|@{\hskip 1em}c@{\hskip 1em}|@{\hskip 1em}c}
h(y)&h(Vy)&h(V^2y)&h(V^3y)\\[.2em]\hline
y_i&y_j&-y_i&-y_j\\[.2em]\hline
y_iy_j&-y_iy_j&y_iy_j&-y_iy_j\\[.2em]\end{array}$$
The rows add up to zero, so
$$\sum_{U\in{\curlyh}}h(Uy)=\frac14\sum_{n=0}^3\sum_{U\in\curlyh}h(V^nUy)=0$$
and the cubic average is zero.
For $1$, the cubic average is $1$.
For $y_i^2$, the cubic average is $(y_1^2 + \cdots + y_d^2)/d$.
So for any monomial of degree two or less,
$$\cuaver h(y) = h(0) + \nabla^2h(0) {|y|^2\over2d}.$$
Add this up for every monomial in $P$ to get the identity
$$\cuaver{P}(y) = P(0) + \nabla^2P(0) {|y|^2 \over 2d}.$$This is what we wanted.\qed

\medskip
\angry Remark. We use the cubic symmetry of the sum in only two places: here and Lemma \ref{diameter}. See the footnote in Section \ref{introduction} about using a different symmetry group.

\end{subsubsection}
\end{subsection}
\begin{subsection}{The first two moves: the bookkeeping}
\label{firsttwomoves}
We can only do the cookie-cutter move when the table set is bulky and has a boundary of measure zero, and the hand sets are small balls. We'll use the first and second moves to get into that situation.
This is possible by the lemmas:

\begin{lemma} Let $\eta > 0$. If $A$ is a bounded open set,
there is a bulky open set $A_0 \subseteq A$ with $\lambda(A \setminus A_0) < \eta$ so that the boundary of $A_0$ has measure zero.\label{secondmoveuse}
\end{lemma}

\angry Proof. The map $t\mapsto\lambda(A^{-t})$ is bounded and monotone, so it's continuous almost everywhere.
Let $t_0$ be a point of continuity for this decreasing function with $\lambda(A^{-t_0}) > \lambda(A) - \eta$.
Then $\lambda(\bigcap_{s < t_0} A^{-s}) = A^{-t_0}$.

Let $A_0$ be the bulky set in the equivalence class of $A^{-t_0}$.
Then $\overline{A_0} \subseteq A^{-s}$ when $s < t_0$, so $\partial A_0 \subseteq (\bigcap_{s < t_0} A^{-s}) \setminus A_0$,
which is the difference of two sets with the same measure.
Therefore, $\lambda(\partial A_0) = 0$.  \qed

\begin{lemma}[Special case of the Vitali covering theorem] Let $B$ be a bounded open set. Let $\eta > 0, R > 0$. There are disjoint open balls $b_1,\ldots,b_m \subseteq B$ with radius less than $R$ so that the measure of $\lambda(B \setminus (b_1 \cup \cdots \cup b_m)) < \eta$.
\label{firstmoveused}
\end{lemma}

\angry Proof. This is well known. See for example \cite{evansgariepy}, Theorem 1.26.\qed

\medskip
Before the $n$-th cookie-cutter move, we'll shrink the table set by a small amount to be chosen later using Lemma~\ref{secondmoveuse}.

When we have to break down the hand into balls smaller than $R$, we'll use the first move and the lemma above to replace all the hand sets with balls. The $n$-th time we do the first move, we choose $\eta = \varepsilon / 2^{n+1}$ in Lemma~\ref{firstmoveused}, so that the total lost mass from the first move is less than $\varepsilon / 2$.

Note that the number of balls in the hand may become very large.

\end{subsection}
\begin{subsection}{The cookie-cutter move always makes progress}
If $E$ is some open set in $\mathbb R^d$, we say that its {\it second moment} is $\int_E |y|^2 \,dy$. This is the same as its `moment of inertia' in two dimensions.

If $A$ is the table set and $B_1, \ldots, B_m$ are the hand sets,
then the \emph{total second moment} is
$\int_A |y|^2 \,dy + \sum_{j=1}^m \int_{B_j} |y|^2 \,dy$.
All the sets in the game are contained in the starting sum $A \oplus B$,
so the total second moment is never more than $(\lambda(A)+\lambda(B))\mathop{\text{rad}}(A \oplus B)^2$.
Here $\rad(E)$ is the radius $\{|x|: x \in E\}$.

We remember some facts about the second moment.
If a ball of radius $r$ is centered at zero, then its second moment is $$\int_{B_r} |y|^2 \,dy = d\lambda(B_1) \int_0^r \rho^2 \times \rho^{d-1}\,d\rho = {d \over d+2} r^2 \lambda(B_r).$$
If the center of mass of a set $E$ is $x$, then its second moment is $$|x|^2 \lambda(E) + \int_{E-x} |y|^2 \,dy.$$

\begin{subsubsection}{The effect of a cookie-cutter move}

The next lemma says essentially that a cookie-cutter move either increases the second moment, or it moves measure outside of the table set.

\begin{lemma}Let $0 < R < \delta/2$. Suppose that we do a cookie-cutter move, smashing a ball $B_r(x)$ with $r < R$ into $\mathscr C(x, R, A)$
to get a new set $E$.

Let $\delta \sigma$ be the change in second moment during the move. Let $\mu = \lambda(B_r)$ be the mass of the ball, and let $\nu = \lambda(E \setminus A)$ be the mass that's moved outside of the set by the cookie-cutter move. Then
$$\delta\sigma + |\curlyh|R^2 \nu \ge {2 \over d+2} R^2\mu.$$\label{odometerlemma}\vskip -1em\vskip 0em\end{lemma}

\angry Proof. The second moment of the ball was
$$\left({d \over d+2}r^2 + |x|^2\right) \mu.$$
The measure of the new set $E$ is the same as the mass of the ball $\mu$, and its centre of mass is $x$ by the cubic symmetry. The second moment of $E - x$ is at least $R^2\lambda((E-x)\setminus{B_R})$, so
$$\int_E |y|^2 \,dy = |x|^2 \mu + \int_{E - x} |y|^2 \,dy \ge |x|^2 \mu + R^2\lambda(E \setminus B_R(x)).$$

So the change in total second moment is
\begin{align*}
\delta \sigma &\ge R^2 \lambda(E \setminus B_R(x)) - {d \over d+2} r^2 \lambda(B_r) \\&\ge R^2\left(\lambda(E \setminus B_R(x)) - {d \over d+2} \lambda(B_r)\right).
\end{align*}

By definition, the set $E$ is disjoint from the cookie-cutter set,
which is $B_R(x)\cap\bigcap_{U\in \curlyh} U_xA$.
Therefore $E \cap B_R(x) \subseteq \bigcup_{U \in \curlyh} (U_xA)^c$,
and $$\lambda(E \cap B_R(x)) \le \lambda\left(E \cap \bigcup_{U \in \curlyh} (U_x A)^c\right) \le \sum_{U \in \curlyh} \lambda(E \setminus U_xA) = |\curlyh| \lambda(E \setminus A).$$
So $\lambda(E \setminus B_R(x)) \ge \lambda(E) - |\curlyh| \lambda(E \setminus A) = \lambda(B_r) - |\curlyh| \lambda(E \setminus A)$.

Substituting this in the inequality above, we get
$$\delta\sigma \ge R^2\left({2 \over d+2} \lambda(B_r) - |\curlyh| \lambda(E \setminus A)\right),$$ and rearranging gives us the result.
\qed

\begin{corollary}Let $R > 0$. Suppose that the mass in the hand is $m$, and every set in the hand is a ball of radius less than $R$. For each ball currently in the hand, we carry out the following steps: \label{winningthesmashgame}

\begin{itemize}
\item[1.] Use the second move to shrink the table set by a small amount so that it's bulky and its boundary has measure zero, as in Section \ref{firsttwomoves}.
\item[2.] Smash $\mathscr{C}(x, R, A)$ into the ball, replacing it by a new set $E$.
\item[3.] Use the fourth move to put $E \setminus \overline{A}$ on the table, leaving $E \cap A$ in the hand.
\end{itemize}

Let the total change in second moment from the cookie-cutter move be $\Delta\sigma$, the total change from the second move be $\Delta \sigma'$, and the total decrease in mass in hand from the fourth move be $\Delta m$. Then $$\Delta \sigma + \Delta \sigma' + |\curlyh| R^2 \Delta m \ge {1 \over d+2} R^2 m.$$The total $s$ integral increases by at most $C_sR^3 m$.\label{smashgame}\end{corollary}

\medskip
\angry Proof. Apply Lemma~\ref{odometerlemma} to each move and add up the inequalities to get
$$\Delta\sigma + |\curlyh| R^2 \Delta m \ge {2 \over d+2} R^2 m.$$
The second moment goes down every time we shrink the table set, but we can make the loss arbitrarily small. If we choose the ``small amount'' in step 1 to be $${1 \over 2^{n+1}}\min\left\{\varepsilon, {R^2m\over(d+2)\mathop{\text{rad}}(A\oplus{B})^2}\right\}$$ where $n$ starts at one at the start of the game and increases every time we use the second move, then $\Delta\sigma' \ge -R^2m/(d+2)$. Adding this inequality to the one above gives the result we wanted.

The total $s$ integral increases by at most $C_s \lambda(B_r) R^3$ at each step, which means that the whole process increases it by at most $C_s m R^3$.\qed

\end{subsubsection}
\end{subsection}
\begin{subsection}{How to win smash game}
\label{endofstrategizing}
We will now give a strategy for smash game which proves by construction that it can always be won.

Recall that we start with a table set $A$, a hand set $B$, a small positive number $\varepsilon$, and a smooth, nonnegative superharmonic function $s$ defined on some set $(A \oplus B)^\delta$, where $\delta > 0$. We have to get the mass in hand below $\varepsilon$ without increasing the total $s$ integral by more than $\varepsilon$, and without decreasing the current mass by more than $\varepsilon$.

\begin{subsubsection}{The strategy}
Let $R_n \in (0, \delta/2N)$ be a sequence satisfying $\sum R_n^2 = \infty$, but $\sum R_n^3 < \varepsilon / 2C_s\lambda(B)$,
where $C_s$ is the constant in Lemma~\ref{x}.
For example, $\sum n^{-3/2}$ is less than $3$, so we could take $R_n$ to be either $\delta / 2N$ or $\varepsilon/6C_s\lambda(B)n^{1/2}$, whichever is smaller.

We repeat the following steps until the mass in hand is below $\varepsilon$. On the $n$-th round:
\begin{itemize}
\item Break each hand set into balls of radius $< R_n$, as in Section \ref{firsttwomoves}.
\item Then carry out the steps in the statement of Corollary \ref{winningthesmashgame} to smash all the balls into cubically symmetric subsets of the table set.
\end{itemize}

The mass in hand at the start of the round is at most $\lambda(B)$, so each round increases the total $s$ integral by at most $C_s \varepsilon R_n^3 \lambda(B)$.
We've chosen the numbers $R_n$ so that the sum of this over all $n$ is less than $\varepsilon$.

The total decrease in the current mass over all rounds is also less than $\varepsilon$, because the losses from the first two moves are bounded by $\sum \varepsilon/2^{n+1} = \varepsilon/2$
and the other two moves don't lose mass. These two paragraphs together tell us that, if we play this way, we'll never lose. The only way we can fail to win is if the game never ends.

\end{subsubsection}
\begin{subsubsection}{The strategy works}
Here we'll prove that the strategy does always win after a finite number of moves.

\begin{lemma}The strategy above always wins smash game.\end{lemma}

\angry Proof. Let the total second moment at the start of the $n$-th round be $\sigma_n$, and similarly let the mass in hand at the start of the round be $m_n$. By Corollary~\ref{winningthesmashgame},
$$\sigma_{n+1} - \sigma_n + |\curlyh| R_n^2 (m_n - m_{n+1}) \ge {1 \over d+2} R_n^2m_n.$$
If we haven't won by time $M$, then $m_n \ge \varepsilon$ for $1 \le n \le M$. The second moment is bounded by
$\sigma_b:=(\lambda(A)+\lambda(B))\mathop{\text{rad}}(A\oplus{B})^2,$ and $R_n < \delta$ and $m_n$ are decreasing with $m_1 = \lambda(B)$, so
$$\sigma_b + |\curlyh| \delta \lambda(B) \ge {\varepsilon \over d+2} \sum_{n=1}^N R_n^2.$$

Remember that we chose the radii $R_n$ so $\sum R_n^2 = \infty$. Let $M$ be so large that $\sum_{n=1}^M R_n^2$ is greater than $(d+2)(\sigma_b + |\curlyh| \delta \lambda(B)) / \varepsilon$.
If the game continues for $M$ rounds, then the above inequality will be violated, which is a contradiction. We never lose the game with our strategy, so the game must have been won before then.
\qed

\end{subsubsection}
\begin{subsubsection}{The sum is a quadrature domain}

We now know that we can win smash game, so we can use Theorem \ref{canwinit} to prove the quadrature domain inequality for smooth superharmonic functions.
\begin{corollary}\label{quad}
Let $A, B$ be bounded open sets. Then $$\int_{A \oplus B} s \,dx \le \int_A s \,dx + \int_B s \,dx$$ for any smooth superharmonic function $s$ defined on a neighbourhood of $\overline{A \oplus B}$.
\end{corollary}

\angry Proof. If $\delta$ is small enough, then the domain of $s$ contains $(A \oplus B)^{2\delta}$. Of course, it's bounded below on any compact subset of its domain, so $c := \min\{s(x): x \in \overline{(A \oplus B)^\delta}\}$ is finite and $s - c$ is a smooth nonnegative superharmonic function on $(A \oplus B)^\delta$.

Start smash game with $A$ on the table, $B$ in the hand, and $s - c$ as the function. Using the strategy above, we can win the game, so Theorem \ref{canwinit} applies and we get $$\int_{A \oplus B} s - c_- \,dx \le \int_A s - c \,dx + \int_B s - c \,dx.$$ By conservation of mass, $\int_{A \oplus B} c\,dx = c \lambda(A \oplus B)$ is the same as $\int_A c \,dx + \int_B c \,dx$, so that part cancels out and we have the inequality that we want.
\qed

\medskip
We want more than that, though: we want the quadrature inequality to hold for all integrable superharmonic functions on $A \oplus B$. However, this follows easily using standard approximation results together with the monotonicity of the sum.

First, we prove the statement for integrable superharmonic functions on a slightly larger domain:
\begin{corollary}\label{superharmoniconnbhd}
The same inequality holds if $s$ is any integrable superharmonic function defined on a neighbourhood of $\overline{A \oplus B}$.
\end{corollary}

\angry Proof. Pick a small number $\delta$ so that the set $(A \oplus B)^\delta$ is contained in the domain of definition of $s$.
Let $C := (A \oplus B)^{\delta/2}$.

Let $\psi$ be a smooth nonnegative bump function which is zero outside the ball of radius one, and let~$s_m = s * [m^d \psi(x/m)]$ for $m > 4/\delta$.
This function is defined for any point in $(A \oplus B)^{\delta/4}$, and on that set, $s_m$ is smooth and superharmonic, as well as nonnegative.
Therefore, $$\int_{A \oplus B} s_m\,dx \le \int_A s_m\,dx + \int_B s_m\,dx$$ by the previous corollary. It's a standard result that $s_m \to s$ in $L^1(A \oplus B)$, so
$$\int_{A \oplus B} s\,dx\le\int_As\,dx+\int_Bs\,dx.$$
This is the inequality that we want to prove.\qed

\bigskip
Now we get the inequality for any integrable superharmonic function.
\begin{theorem}\label{superharmonicultimate}
Let $A, B$ be bounded open sets. Then $\int_{A \oplus B} s \,dx \le \int_A s \,dx + \int_B s\,dx$ for any integrable superharmonic function $s$ on $A \oplus B$, or in other words, $A \oplus B$ is a quadrature domain for $\1_A+\1_B$.
\end{theorem}

\angry Proof. We use the above corollary to get the inequality on smaller sets and then use dominated convergence as the sets increase. If $x \in A^{-\varepsilon}$, then $x$ is farther than $\varepsilon$ from the boundary, so $(A^{-\varepsilon})^{\varepsilon} \subseteq A$. We recall the inflation inclusion $(\ref{inflationinclusion})$, which says that $(A \oplus B)^\varepsilon \subseteq A^\varepsilon \oplus B^\varepsilon$. By monotonicity,
$$(A^{-\varepsilon} \oplus B^{-\varepsilon})^\varepsilon \subseteq (A^{-\varepsilon})^\varepsilon \oplus (B^{-\varepsilon})^\varepsilon \subseteq A \oplus B.$$
Set $C_\varepsilon := A^{-\varepsilon} \oplus B^{-\varepsilon}$. Then the inclusion above tells us that $(C_\varepsilon)^\varepsilon \subseteq A \oplus B$, which means that $A \oplus B$ contains a neighbourhood of $\overline{C_\varepsilon}$. By Corollary \ref{superharmoniconnbhd},
$$\int_{C_\varepsilon} s\,dx \le \int_{A^{-\varepsilon}} s\,dx + \int_{B^{-\varepsilon}} s\,dx.$$
Let $C := \bigcup_n C_{1/n}$. Then $s \1_{C_{1/n}} \to s \1_C$ pointwise, and similarly $s \1_{A^{-\varepsilon}} \to s \1_A$ and $s \1_{B^{-\varepsilon}} \to s \1_B$. All the functions are dominated by $|s|$, which by assumption is integrable on $A \oplus B$. Therefore, by dominated convergence,
$$\int_{C} s\,dx \le \int_A s\,dx + \int_B s\,dx.$$
By conservation of mass, $\lambda(C_\varepsilon) = \lambda(A^{-\varepsilon}) + \lambda(B^{-\varepsilon}) \to \lambda(A) + \lambda(B)$ as $\varepsilon \to 0$, so $\lambda(C) = \lambda(A) + \lambda(B)$, and $C$ is contained in $A \oplus B$.
Therefore, they are essentially equal, and $\int_C s \,dx = \int_{A \oplus B} s \,dx$. So $$\int_{A\oplus{B}}s\,dx\le\int_As\,dx+\int_Bs\,dx = \int (\1_A+\1_B) \,s\,dx$$for any integrable superharmonic function on $A\oplus{B}$, which is what we have claimed.

We recall that this is the definition of a quadrature domain for $\1_A+\1_B$: $$\int_{A \oplus B} s\,dx \le \int_A s \,dx + \int_B s \,d = \int s w \,dx$$ where $w := \1_A + \1_B$.\qed

\end{subsubsection}
\end{subsection}
\begin{subsection}{Conclusion}

\begin{subsubsection}{There's no other sum of open sets}

\begin{theorem}The Diaconis-Fulton smash sum is the only sum of open sets that satisfies the six requirements plus bulkiness.\end{theorem}

\angry Proof. Let $\oplus$ be any sum of open sets satisfying the requirements, and $A$ and $B$ any two bounded open sets.
By Theorem~\ref{superharmonicultimate}, $A \oplus B$ is a quadrature domain for $\1_A + \1_B$.

By Theorem~\ref{superharmonicultimate}, $A \oplus B$ is a quadrature domain for $\1_A + \1_B$, and so is the Diaconis-Fulton smash sum of $A$ and $B$. Quadrature domains are essentially unique by Theorem~\ref{sakaisx}, so the two sets are essentially equal, and they are both bulky, so they are really equal by Lemma~\ref{really}.

\begin{corollary}If a sum of open sets $\oplus$ satisfies the six requirements in the essential sense, then for any bounded open sets $A,B$, the sum of $A$ and $B$ is essentially equal to the Diaconis-Fulton smash sum of $A$ and $B$.\end{corollary}

\angry Proof. By Lemma~\ref{aseventhaxiom}, the bulked sum $(A, B) \mapsto \bulk{A \oplus B}$ satisfies the six requirements in the strong sense, plus bulkiness. Therefore, by the theorem, the bulked sum is the Diaconis-Fulton smash sum, and $A \oplus B$ is essentially equal to $\bulk{A \oplus B}$.
\qed

\end{subsubsection}
\begin{subsubsection}{Some open questions}

Are there sums that satisfy the six requirements in the strong sense, but don't satisfy the requirements of bulkiness?
In particular, is there a sum of open sets with $A \oplus B = A \cup B$ when $A, B$ are disjoint?
It would have to be essentially equal to smash sum, but it's not impossible that sets of measure zero could be left out according to some clever scheme so that the requirements are still satisfied.

\medskip
Let $f(r, s) = (r^d + s^d)^{1/d}$.
We delete the conservation of mass requirement, and instead add:

\begin{itemize}
\item Continuity. If $\lambda(A_n \setdelta A) \to 0$ and $\lambda(B_n \setdelta B) \to 0$, then the measure of the differences $(A_n \oplus B_n) \setdelta (A \oplus B)$ goes to zero.
\item Addition of concentric balls. If $r, s \ge 0$, then the sum $B_r \oplus B_s$ is $B_{f(r, s)}$.
\end{itemize}

It's not hard to show that these are consequences of the six requirements, so this new set is weaker.
Is there still only one sum satisfying them?

\medskip
We can also change the function $f$. For example, if we set $f(r, s) = \max\{r, s\}$, then the union sum $A, B \mapsto A \cup B$ satisfies the above requirements.
Are there any other functions $f$ for which a sum exists?

\medskip
Could one develop a similar uniqueness result for the sum on a general Riemannian manifold, as it appears in \cite{riemannian}? Full translation invariance would be impossible unless the manifold had constant curvature, but one could ask for the sum to be approximately symmetric for small sets.
\end{subsubsection}
\end{subsection}
\end{section}

\bibliographystyle{acm}
\bibliography{smash-sum}
\end{document}